\documentclass[review]{elsarticle}
\newtheorem{Theorem}{Theorem}[section]

\newtheorem{Proposition}[Theorem]{Proposition}

\newtheorem{Definition}[Theorem]{Definition}
\newtheorem{Example}[Theorem]{Example}

\newtheorem{Remark}[Theorem]{Remark}

\newproof{pf}{Proof}

\usepackage{amsmath}
\usepackage{amssymb}
\usepackage{amsfonts}

\usepackage{lineno,hyperref}
\usepackage{color}

\newcommand{\precN}{\ensuremath{\stackrel{\mathrm{N}}{\prec}}}

\modulolinenumbers[5]
\journal{Journal of \LaTeX\ Templates}
\begin{document}
\begin{frontmatter}
\title{Genera and minors of multibranched surfaces}
\author[mymainaddress]{Shosaku Matsuzaki}
\ead{shosaku@aoni.waseda.jp}
\author[mysecondaryaddress]{Makoto Ozawa}
\ead{w3c@komazawa-u.ac.jp}
\address[mymainaddress]{Integrated Arts and Sciences, Waseda University,\\ 1-6-1 Nishiwaseda, Shinjuku-ku, Tokyo, 169-8050, Japan}
\address[mysecondaryaddress]{Department of Natural Sciences, Faculty of Arts and Sciences, Komazawa University,\\ 1-23-1 Komazawa, Setagaya-ku, Tokyo, 154-8525, Japan}

\begin{abstract}
We say that a $2$-dimensional CW complex is a {\em multibranched surface} if we remove all points whose open neighborhoods are homeomorphic to the $2$-dimensional Euclidean space $\mathbb{R}^2$, then we obtain a $1$-dimensional complex which is homeomorphic to a disjoint union of some $S^1$'s. We define the genus of a multibranched surface $X$ as the minimum number of genera of $3$-dimensional manifold into which $X$ can be embedded. We prove some inequalities which give upper bounds for the genus of a multibranched surface.
A multibranched surface is a generalization of graphs. Therefore, we can define ``minors'' of multibranched surfaces analogously. We study various properties of the minors of multibranched surfaces.
\end{abstract}

\end{frontmatter}

\section{Introduction}\label{intro}
It is a fundamental problem to determine whether or not there exists an embedding from a topological space into another one.
The Menger--N\"{o}beling theorem ({\cite[Theorem 1.11.4.]{E}}) shows that any finite $2$-dimensional CW complex can be embedded into the $5$-dimensional Euclidian space $\mathbb{R}^{5}$.
This is a best possible result since for example, the union of all $2$-faces of a $6$-simplex cannot be embedded in $\mathbb{R}^{4}$ ({\cite[1.11.F]{E})}.
But, if the subspace of a finite $2$-dimensional CW complex consisting of all points which do not have an open neighborhood homeomorphic to $\mathbb{R}^{2}$ is a possibly disconnected $1$-dimensional manifold, then the CW complex can be embedded into $\mathbb{R}^{4}$ (Proposition \ref{prop_r4}).
We call such CW complexes {\em multibranched surfaces}, and moreover obtain a necessary and sufficient condition for a multibranched surface to be embeddable into some orientable closed $3$-dimensional manifold (Proposition \ref{prop_embedding}).
This is a starting point to study such multibranched surfaces via $3$-dimensional manifolds.

In this paper, we introduce the {\em genus} of such a multibranched surface as the minimal Heegaard genus of $3$-dimensional manifolds into which the multibranched surface can be embedded.
In section \ref{sec_genus}, we will give some upper bounds for the genus of a multibranched surface (Theorem \ref{thm_upp_of_mingenus}, \ref{thm_min_genus}).
In section \ref{sec_homology}, we will describe the first homology groups of multibranched surfaces, and calculate for some examples.
It can be used to determine whether or not a multibranched surface can be embedded into the $3$-sphere.
This constructively explains more details of the calculation by using the determinant of some matrix in \cite{EMO}.
In \cite{EMO}, we also studied the criticality of a multibranched surface for the $3$-sphere $S^3$ and have given some critical multibranched surfaces for $S^3$.
As the Kuratowski's theorem (\cite{K}) characterized the $2$-sphere $S^2$ by means of the obstruction set of graphs, it might be possible to characterize a closed $3$-dimensional manifold by means of some obstruction set of multibranched surfaces.

In Graph Theory, Robertson and Seymour introduced the minor theory which gives a most important structure on the set of graphs.
Since we can regard a graph as a $1$-dimensional multibranched manifold,
it would be natural to consider a similar minor theory for multibranched surfaces.
In Section 5, we will define the minor for multibranched surfaces and introduce some intrinsic properties which are minor closed.
Thus we arrive at the obstruction set for those intrinsic properties,
and give some examples which belong to the obstruction set.
And also we define the neighborhood minor for multibranched surfaces.
The neighborhood minor sets a preorder on the set of multibranched surfaces, and behaves well on some basic operations
(Proposition \ref{minorclosed}).
In particular, if $X$ is a neighborhood minor of $Y$, then the genus of $X$ is less than or equal to that of $Y$.

To summarize this paper, we have found a well-behaved class of $2$-dimensional CW complexes (which we call {\em regular multibranched surfaces}) and derived an invariant of regular multibranched surfaces from the Heegaard genus of $3$-dimensional manifolds which is known to be the most fundamental invariant of $3$-dimensional manifolds.
In the future, we expect some characterization of each $3$-dimensional manifold by means of the obstruction set of regular multibranched surfaces like as the Kuratowski's theorem.

\section{Preliminaries}\label{preliminaries}

\subsection{Multibranched surfaces}\label{mbs}
\begin{figure}[htbp]\begin{center}\scalebox{0.9}{\includegraphics*{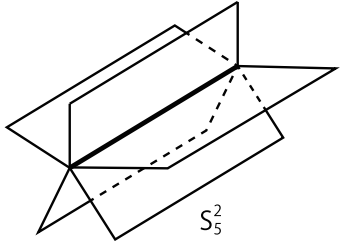}}\end{center}\caption{an open neighborhood homeomorphic to $S^2_5$}\label{S25}\end{figure}

Let $\mathbb{R}^n_+$ be the upper half space $\left\{ (x_1, x_2, \ldots, x_n) \in \mathbb{R}^n  |  x_n \geq 0 \right\}$ of $n$-dimensional Euclidean space $\mathbb{R}^n$.
The quotient space obtained from $i$ copies of $\mathbb{R}^n_{+}$ by identifying with their boundaries, $\partial{\mathbb{R}^n_+}= \left\{ (x_1, x_2, \ldots, x_n)  \in \mathbb{R}^n |  x_n=0 \right\}$, is denoted by $S^{n}_i$. We note that $S^{n}_2$ is homeomorphic to $\mathbb{R}^n$.
See Figure \ref{S25}.

\begin{Definition}\label{multibranchedmfd}{\rm
An $n$-dimensional CW complex $X$ is an $n$-{\it dimensional multibranched manifold}
if for every point $x \in X$ there exist a positive integer $i$ and an open neighborhood $U$ of $x$ such that $U$ is homeomorphic to $S_i^n$.
}\end{Definition}

We call a $2$-dimensional multibranched manifold a {\it multibranched surface}.
In this paper, we consider multibranched surfaces which are constructed by gluing some compact $2$-dimensional manifolds into their boundaries. 

\noindent We prepare a closed $1$-dimensional manifold $L$, a compact $2$-dimensional manifold $E$ and a continuous map $\phi: \partial E \to L$ satisfying the following conditions.
\begin{enumerate}
\item For every connected component $e$ of $E$, $\partial e \not= \emptyset$.
\item For every connected component $c$ of $\partial E$, the restriction $ \phi|_c : c \to \phi(c) $ is a covering map.
\end{enumerate}

\noindent The quotient space $X=L \cup_\phi E$ is called the {\it multibranched surface} obtained from the triple $(L, E; \phi)$.
We note that $L$ and $E$ are not necessarily connected.

A connected component of $L$ (resp. $E$, $\partial E$)  is said to be a {\it branch} (resp. {\it sector}, {\it prebranch}) of $X$.
We note that every branch of $X$ is homeomorphic to the $1$-sphere $S^1$.
The set consisting of all branches (resp. sectors) is denoted by $\mathcal{L}(X)$ (resp. $\mathcal{E}(X)$).

For a prebranch $c$ of a multibranched surface $X$, the covering degree of the covering map $\phi|_c:c\to \phi(c)$ is called the {\it degree} of $c$, denoted by $d(c)$. We note that $d(c)$ is a positive integer.
A prebranch $c$ of $X$ is said to be {\it attached} to a branch $l$ if $\phi(c)=l$.
The number of prebranches which are attached to a branch $l$ of $X$ is called the {\it index} of $l$, denoted by $i(l)$.


\begin{Definition}\label{def_decomp}{\rm
Let $X=L \cup_\phi E$ be a multibranched surface.
By cutting $E$ along $\partial N(\partial E;E)-\partial E$ and pasting
two disks, we obtain disks (denoted by $\hat{E}$) whose boundaries
coincide with that of $E$ and closed surfaces $F$.
Put $\hat{X} = L \cup_\phi \hat{E}$.
Then, $X$ is decomposed into $\hat{X}$ and $F$ and conversely $X$ is
obtained from $\hat{X}$ and $F$ by tubings.
We call this decomposition a {\it standard decomposition} of $X$, and
denote it by $X=\hat{X} \# F$.
}\end{Definition}

\begin{Proposition}\label{prop_r4}
Every multibranched surface is embeddable into the $4$-dimensional Euclidean space $\mathbb{R}^4$.
\end{Proposition}

\begin{figure}[htbp]\begin{center}\scalebox{0.9}{\includegraphics*{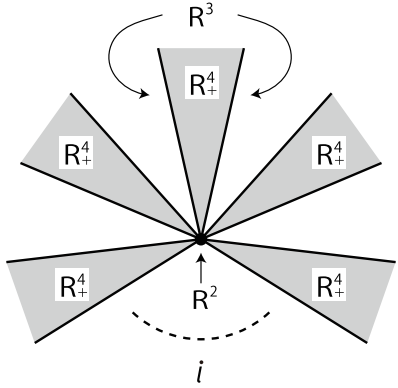}}\end{center}\caption{A simplicial complex $K$ in $\mathbb{R}^4$.}
\label{r4}\end{figure}

\begin{pf}
Let $X$ be a multibranched surface.
Let $K$ be a simplicial complex embedded into $\mathbb{R}^4$
such that
$K$ is obtained by identifying with $\mathbb{R}^2$ part of $i$ copies of $\mathbb{R}^4_+$ with $\mathbb{R}^2  \subset \mathbb{R}^3 = \partial \mathbb{R}^4_+ \subset \mathbb{R}^4_+$, where $i={\rm max}\{i(l)| l \mbox{ is a branch of } $X$\}$.
See in  Figure \ref{r4}.
We consider a standard decomposition $X=\hat{X} \# F$ of $X$.\\
(Step 1.)
First, we embed $\hat{X}$ into $K$ as follows.
Every branch of $\hat{X}$ is contained in the part $\mathbb{R}^2$ of $K$.
Every component of $\hat{X} -L $ is embedded into the each part $\mathbb{R}^4_+ - \mathbb{R}^2$ of $K$
and every component $c$ of $\partial \left( N(\partial E; E) \right)-L$ is trivial knot contained in $\mathbb{R}^3 = \partial \mathbb{R}^4_+ \subset \mathbb{R}^4_+$.
Then, we obtain an embedding of the multibranched surface $\hat{X}$ into $K$.\\
(Step 2.)
Second, we embed $F$ into $\mathbb{R}^4$ with $F \cap K = \emptyset$.
By tubing $F$ and $\hat{X} \subset K$, we obtain an embedding of $X=\hat{X} \# F$ into $\mathbb{R}^4$.
$\Box$
\end{pf}

We give an orientation for each branch and each prebranch $c$ of $X$.
The {\it oriented degree} of a prebranch $c$ of $X$ is defined as follows:
if the covering map $\phi|_c:c\to \phi(c)$ is orientation preserving,
the {\it oriented degree} $od(c)$ of $c$ is defined by $od(c)=d(c)$ and
if it is orientation reversing, the oriented degree is defined by $od(c)=-d(c)$.

\begin{Definition}\label{regular}{\rm
A multibranched surface $X$ is {\it regular} if the following condition is satisfied.
\begin{enumerate}
\item [] For every branch $l$ and every prebranch $c$ and $c'$ of $X$ which are attached to $l$, $d(c)=d(c')$.
\end{enumerate}
}\end{Definition}

Let $X$ be a regular multibranched surface.
Since each pair of prebranches $c, c'$ of $X$ which are attached to a branch $l$ has same degree,
then we define the {\it degree} of a branch $l$ as $d(l)=d(c)=d(c')$.

In this paper,
the cardinality of a set $S$ is denoted by $\# S$.
For a $3$-dimensional manifold $M$, we will denote by $\mathring{M}$ the set $M - \partial M$.

\subsection{Neighborhoods of multibranched surfaces}\label{neighborhood}

For a regular multibranched surface $X$, we define the ``circular permutation system'' and ``slope system'' of $X$ as follows.
A circular permutation of prebranches which are attached to a branch $l$ is called a {\it circular permutation} of $l$.
A collection $\mathcal{P}=\{ \mathcal{P}_l \}_{l \in \mathcal{L}(X)}$ is called a {\it circular permutation system} of $X$ if $\mathcal{P}_l$ is a circular permutation of $l$.
For a branch $l$, a rational number $p/q$ with $q=d(l)$ is called a {\it slope} of $l$.
A collection $\{ \mathcal{S}_l \}_{l \in \mathcal{L}(X)}$ is called a {\it slope system} of $X$ if $\mathcal{S}_l$ is a slope of $l$.

We define a compact $3$-dimensional manifold with boundary, called a ``neighborhood'' of a regular multibranched surface $X$ in the Definition \ref{nbh}.
This is uniquely determined up to homeomorphism by a pair of a circular permutation system $ \mathcal{P} = \{ \mathcal{P}_l \}_{ l \in \mathcal{L}(X) } $ and a slope system $\mathcal{S}=\{ \mathcal{S}_l \}_{ l \in \mathcal{L}(X) }$ of $X$.

\begin{Definition}\label{nbh}{\rm
Let $X$ be a regular multibranched surface and let $\mathcal{P}=\{ \mathcal{P}_l \}_{l \in \mathcal{L}(X)}$ and $\mathcal{S}=\{ \mathcal{S}_l \}_{ l \in \mathcal{L}(X) }$ be a permutation system and a slope system of $X$ respectively.
We will construct the $3$-dimensional manifold by the following procedure.
First, for each branch $l$ of $X$ and each sector $e$ of $X$ we take a solid torus $l \times D^2$, where $D^2$ a disk and take a product $e \times [-1, 1]$. (If $e$ is non orientable, we take a twisted $I$-bundle $e \tilde{\times} [-1, 1]$ over $e$.)
We give orientations for these $3$-dimensional manifolds.
Next, we glue them depending on the permutation system $\mathcal{P}$ and the slope system $\mathcal{S}$
by assigning of the slope $\mathcal{S}_l$ of $l$ to the isotopy class of a loop $k$ in $\partial ( l \times D^2)$
by an orientation reversing homeomorphism $\Phi: c \times [-1, 1] \to N \left(k; \partial \left(  l \times D^2 \right) \right)$ or $\tilde{\Phi}: c \tilde{\times} [-1, 1] \to N \left(k; \partial \left(  l \times D^2 \right) \right)$ with $\phi(c)=l$, $c \subset \partial e$.
See Figure \ref{neighborhood}.
Then, we uniquely obtain a compact and orientable $3$-dimensional manifold with boundary,
denoted by $N(X; \mathcal{P}, \mathcal{S})$.
The $3$-dimensional manifold $N(X; \mathcal{P}, \mathcal{S})$ is called the {\it neighborhood} of $X$ with respect to $\mathcal{P}$ and $\mathcal{S}$. 
We note that $X$ has a neighborhood if and only if $X$ is regular.
The set consisting of all neighborhoods of $X$ is denoted by $\mathcal{N}(X)$.
}\end{Definition}

\begin{figure}[htbp]
	\begin{center}
	\includegraphics[trim=0mm 0mm 0mm 0mm, width=.8\linewidth]{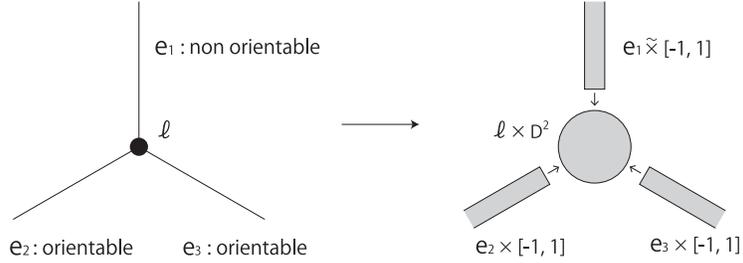}
	\end{center}
	\caption{Neighborhood}
	\label{neighborhood}
\end{figure}

We obtain the following proposition
which assures us of studying multibranched surfaces via invariants of $3$-dimensional manifolds.

\begin{Proposition}\label{prop_embedding}
A multibranched surface is embeddable in some orientable closed $3$-dimensional manifold if and only if the multibranched surface is regular. 
\end{Proposition}

\begin{pf}
($\Leftarrow$)
If a multibranched surface is regular,
there exists a neighborhood $N \in \mathcal{N}(X)$ of $X$.
By attaching some handlebodies to the boundary of $N$, we obtain an orientable closed $3$-dimensional manifold.\\
($\Rightarrow$)
If $X$ is embeddable into an oriented closed $3$-dimensional manifold $M$, there exist a circular permutation system $\mathcal{P}$ and a slope system $\mathcal{S}$ of $X$ such that $N(X; \mathcal{P}, \mathcal{S})$ is homeomorphic to $N(X;M)$.
$\Box$
\end{pf}

\begin{Definition}\label{a_d_graph}{\rm
Let $X$ be a regular multibranched surface.
Suppose that we have fixed a circular permutation system $\mathcal{P}$ of $X$.
For every slope system $\mathcal{S}$ of $X$,
we can glue handlebodies on the boundary  of the neighborhood $N(X; \mathcal{P}, \mathcal{S})$ of $X$ and obtain the closed $3$-dimensional manifold $M$.
We consider the geometric dual graph $G$ of $N(X; \mathcal{P}, \mathcal{S})$ in $M$.
Since the graph $G$ is depending on a circular permutation system $\mathcal{P}$ and  independent of all slope systems of $X$,
then we can define the {\it abstract dual graph} $G(X; \mathcal{P})$ of $X$ with respect to $\mathcal{P}$ by the graph $G$.
The set consisting of all abstract dual graphs of a regular multibranched surface $X$ is denoted by $\mathcal{G}(X)$. 
}\end{Definition}


\section{Genera of multibranched surfaces}\label{sec_genus} 

\subsection{Minimum genera and maximum genera of multibranched surfaces}\label{subsec_genus}
For an orientable $3$-dimensional manifold $M$ with boundary, the minimum Heegaard genus of orientable closed $3$-dimensional manifolds into which $M$ is embeddable is denoted by $eg(M)$, called the {\it embeddable genus} of $M$.

\begin{Remark}\label{rem_genus}{\rm
By the proof of Proposition \ref{prop_embedding}, we have $eg(M) \leq g(M)$, where $g(M)$ is the Heegaard genus of $M$.
}\end{Remark}

\begin{Definition}\label{def_genus}{\rm
For a regular multibranched surface $X$, we define the {\it minimum genus} $g(X)$ and {\it maximum genus} $G(X)$ respectively as following.
\begin{enumerate}
\item [] $g(X)={\rm min} \{ eg(N)  |  N \in \mathcal{N}(X) \}$, \hspace{4mm} $G(X)={\rm max} \{ eg(N) | N \in \mathcal{N}(X) \}$.
\end{enumerate}
}\end{Definition}

We find similar definitions in Graph Theory (for example \cite{BR}, \cite{FGM}).
The next two proposition show that this invariants of regular
multibranched surfaces are non-trivial, and that the gap between them
can be arbitrary large.

\begin{Proposition}\label{prop_gxn}
For every positive integer $n$, there exists a regular multibranched surface $X$ such that $g(X)=G(X)=n$.
\end{Proposition}

\begin{pf}
Let $\bar{X}_n$ be a regular multibranched surface $\bar{X}(p_1,p_2,\ldots,p_n)$ shown in Figure \ref{seifert}, where $l_1$, $l_2$, $l_3$, $\ldots$, $l_n$ are branches of $\bar{X}_n$ and $D_1$, $A_1$, $A_2$, $A_3$, $\ldots$, $A_{n-1}$ are sectors of $\bar{X}_n$ and $|p_i|$ is the degree of $l_i$ and $|p_i|\ge 2$ $(i=1,2,\ldots,n)$. We will show that $g(\bar{X}_n)=G(\bar{X}_n)=n$.

\begin{figure}[htbp]\begin{center}\includegraphics[trim=0mm 0mm 0mm 0mm, width=.9\linewidth]
{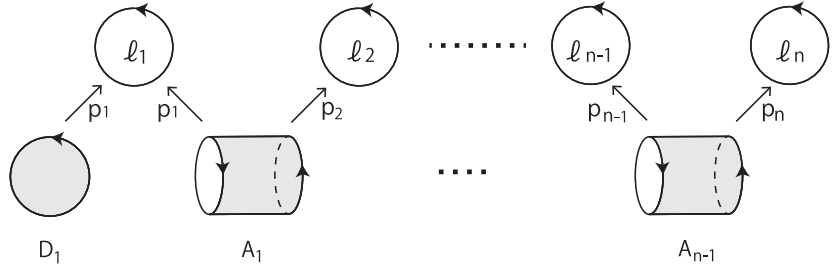}\end{center}
\caption{$\bar{X}(p_1,p_2,\ldots,p_n)$}
\label{seifert}
\end{figure}

\noindent
Suppose that $\bar{X}_n$ can be embedded into a closed orientable $3$-dimensional manifold $M$.
Since $\partial N(\bar{X}_n)$ is homeomorphic to the $2$-sphere
	and $N(\bar{X}_n)=N(\bar{X}_n; \mathcal{P}, \mathcal{S}) = \#_{i=1}^n L(p_i,*)-\mathring{B^3}$
	for every circular permutation system $\mathcal{P}$ and slope system $\mathcal{S}$ of $\bar{X}_n$,
	we have by the additivity of Heegaard genus (\cite{H}) that $g(M)\ge g(N(\bar{X}_n)\cup B^3)=g(\#_{i=1}^n L(p_i,*))=n$.
$\Box$\end{pf}

\begin{Proposition}\label{prop_gxgx}
For positive integer $n$,
there exists a regular multibranched surface $X$ such that $G(X)-g(X) > n$.
\end{Proposition}

\begin{pf}
Let $\Gamma$ be a rose with $2n$ petals $(n\ge 1)$.
See Figure \ref{roses}.
Take a product $\Gamma \times S^1$ and glue a disk $D$ with (vertex of $\Gamma$)$\times S^1$ along its boundary.
We consider this multibranched surface $X=(\Gamma \times S^1)\cup D$.
Let $N_i(X)$ be a neighborhood of $X$ which is determined by a circular permutation of $\Gamma_i$ $(i=1,2)$.

\begin{figure}[htbp]\begin{center}\includegraphics[trim=0mm 0mm 0mm 0mm, width=.5\linewidth]
{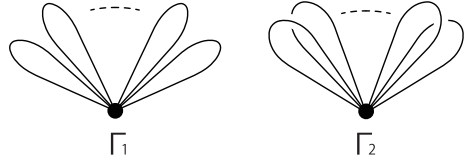}\end{center}
\caption{Roses $\Gamma_1$ and $\Gamma_2$.}
\label{roses}
\end{figure}

Since $\Gamma_1$ is planar with respect to the circular permutation, $\Gamma\times S^1$ can be embedded in the standard solid torus in $S^3$ so that (vertex of $\Gamma$)$\times S^1$ bounds a disk $D$ in the outside of the solid torus.
Thus we have an embedding of $N_1(X)$ into $S^3$ and we have $g(N_1(X))=0$.
On the other hand, we note that $N_2(X)$ is homeomorphic to a product manifold of once punctured orientable surface of genus $n$, say $F$, and $S^1$ to which $2$-handle $N(D)$ attached.
Therefore we have that $\pi_1(N_2(X))\cong \pi_1(F)*\mathbb{Z}/\langle \partial D \rangle \cong \pi_1(F)$
and that $g(N_2(X))\ge rank(\pi_1(N_2(X)))=2n$.
$\Box$
\end{pf}

The next theorem gives a fundamental inequality between the genus of a
multibranched surface and the numbers of its branches and sectors.

\begin{Theorem}\label{thm_upp_of_mingenus}
Let $X$ be a regular multibranched surface. Then,
\[
g(X) \leq \#\mathcal{L}(X) + \#\mathcal{E}(X).
\]
\end{Theorem}

\begin{pf}
Let $X'$ be a regular multibranched surface which is obtained from $X$ by reducing the degree at each branch into $1$.
We will show in the following three steps that $X$ can be embedded into a closed orientable $3$-dimensional manifold of genus $\#\mathcal{L}(X)+ \#\mathcal{E}(X)$.\\
(Step 1: $X'\looparrowright S^3$). We construct an immersion of $X'$ into $S^3$ such that $l_1\cup\cdots \cup l_m$ is the trivial link, where $m=\# \mathcal{L}(X)$, and all multi points are double points contained in the interior of sectors as follows. 
First we decompose $\mathcal{E}(X')$ into a collection $Y=\{D_1,\ldots,D_p\}$ of disks whose boundaries are $l_1\cup\cdots \cup l_m$ and a collection $Z=\{F_1,\ldots,F_q\}$ of closed orientable surfaces, where $\mathcal{E}(X)$ can be obtained by tubing $Y$ and $Z$. 
Then we obtain a standard decomposition $X'=Y\# Z$ of $X'$.
See Figure \ref{fig_YandZ}.
Next we embed $Y$ into $S^3$ so that $l_1\cup\cdots \cup l_m$ is the trivial link, and embed $Z$ into $S^3$ so that $Z\cap Y=\emptyset$.
Since $X'$ can be obtained from $Y$ and $Z$ by tubings, we obtain an immersion of $X'$ into $S^3$ such that all multi points are double points contained in the interior of sectors.\\

\begin{figure}[htbp]\begin{center}
\includegraphics[width=11cm]{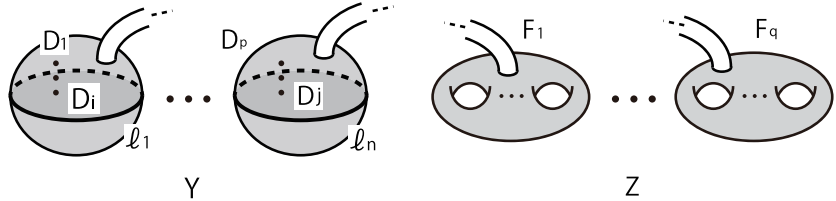}\end{center}
\caption{A standard decomposition $X'=Y\# Z$ of $X'$.}
\label{fig_YandZ}
\end{figure}

\noindent
(Step 2: $X'\hookrightarrow \#^n (S^2\times S^1)$). We construct an embedding of $X'$ into a connected sum of $n$ $S^2\times S^1$ as follows, where $n=\# \mathcal{E}(X)$.
By Step 1, we can take mutually disjoint disks $\Delta_1,\ldots,\Delta_n$ on each sector which contain all double points of $X'$.
Then by performing $0$-slope Dehn surgeries along $\partial \Delta_1,\ldots, \partial \Delta_n$, we obtain a connected sum of $n$ $S^2\times S^1$ into which $X'$ can be embedded.
See Figure \ref{fig_immersion}.

\begin{figure}[htbp]\begin{center}
\includegraphics[width=6cm]{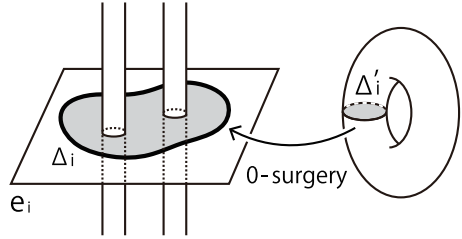}\end{center}
\caption{$0$-slope Dehn surgeries along $\partial \Delta_1,\ldots, \partial \Delta_n$.}
\label{fig_immersion}
\end{figure}

\noindent
(Step 3: $X\hookrightarrow \#^n (S^2\times S^1) \#^m$(Lens space)$=M$). We construct an embedding $X$ into a connected sum of $m$ Lens spaces and $n$ $S^2\times S^1$ as follows.
We perform Dehn surgeries along branches $l_1,\ldots, l_m$ of $X'$ so that $X$ can be recovered.

Indeed, there exists a homeomorphism of a torus which sends $0$-slope
(corresponding to $\partial D_i$) to $p_i/q_i$-slope (the original
slope of a branch $l_i$ in $X$). Therefore, if we remove $N(l_i)$ from
$\#^n (S^2\times S^1)$ and glue a solid torus so that the
homeomorphism is realized, then we have the desired Dehn surgery.
Then $X$ can be embedded into a connected sum of $m$ Lens spaces and $n$ $S^2\times S^1$.
Hence we obtain an embedding $X$ into a closed orientable $3$-dimensional manifold $M$ with $g(M)=m+n$. $\Box$
\end{pf}

In the next theorem,
we show another inequality for the upper bound of the genus of a multibranched surface (cf. Remark \ref{rem_genus}).

\begin{Theorem}\label{thm_min_genus}
Let $X$ be a regular multibranched surface.
Then, for every $N \in \mathcal{N}(X)$,
\[
eg(N) \leq g(\partial N)+\beta (G_N),
\]
where $g(\partial N)$ denotes the sum of genera of the boundary of $N$ and $\beta (G_N)$ is the first betti number of an abstract dual graph of $N$.
Hence we have $g(X) \leq {\rm min} \left\{ g(\partial N) + \beta (G_N) | N \in \mathcal{N}(X) \right\}$.
\end{Theorem}

\begin{pf}
We proceed along the proof in {\cite[p240, 2.\!\! Theorem]{R}}.
As we have seen before, $N$ can be embedded into a closed orientable $3$-dimensional manifold $M$ which is obtained from $N$ by gluing handlebodies $V_1,\ldots,V_{|\partial N|}$ along their boundaries.
For each sector $e_j\in \mathcal{E}$, we take a disk $\delta_j$ in $int e_j$.
Then we obtain a handlebody $V$ which is obtained from $V_1,\ldots,V_{|\partial N|}$ by adding every 1-handle $\delta_j\times [-1,1]$.
We remark that the genus of $V$ is equal to $g(\partial N)+\beta (G_N)$.
The rest $3$-dimensional manifold $W=M-int V$ is also a handlebody since $W$ is obtained from a disjoint union of solid tori $N(l_1 \cup\cdots\cup l_m)$ by adding 1-handles $\bigcup_{i,j} N(\gamma_i^j)$ where $\gamma_i^j$ denotes an arc on a sector $e_j$ such that $e_j-int \delta_j$ is homeomorphic to $N(\partial e_j\cup \bigcup_i \gamma_i^j; e_j)$.
See Figure \ref{fig_handlebody}.
Hence we obtain a Heegaard splitting $M=V\cup W$ of genus $g(\partial N)+\beta (G_N)$.
$\Box$

\begin{figure}[htbp]\begin{center}
\includegraphics[width=12cm]{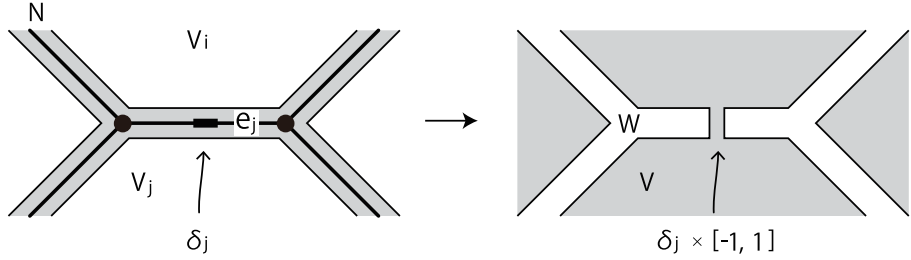}\end{center}
\caption{A handlebody $V$ and $W$.}
\label{fig_handlebody}
\end{figure}
\end{pf}

\section{First homology groups of regular multibranched surfaces}\label{sec_homology}

Give a regular multibranched surface, we can calculate the first homology group of it.
For a branch $l$ and a sector $e$ of $X$,
we define $d(l;e)=\sum_{c \subset \partial e} od(c)$, where $c$ is a prebranch attached to $l$.
(If there exists no prebranches $c$ attached to $l$ with $c \subset \partial e$,
we define $d(l;e)=0$.)

The multibranched surface obtained by the removing a open disk from each sector
is denoted by $\dot{X}$.

\begin{Theorem}\label{thm_h1x}
Let $X$ be a connected and regular multibranched surface. Then,
\[
H_1(X)=
\left[
l_1, \ldots, l_n :
\sum_{k=1}^{n}d(l_k;e_1)l_k, \ldots, \sum_{k=1}^{n}d(l_k;e_m)l_k \right]
\oplus
\mathbb{Z}^{r'(X)}
\],
where $\mathcal{L}(X)=\left\{ l_1, \ldots, l_n \right\}$,
$\mathcal{E}(X)=\left\{ e_1, \ldots, e_m \right\}$ and $r'(X)=r(\dot{X})-n$.
\end{Theorem}

\begin{pf}
Let $O_i$ be an open disk of $X$ which is contained in the sector $e_i$ and let $D_i=\bar{O_i}$ ($1 \leq i \leq m$).
We consider the Mayer-Vietoris sequence of $X_i=\dot{X} \cup D_i$.
\begin{enumerate}
\item []
$\cdots
\overset{\partial_*}{\longrightarrow} H_1 ( \dot{X} \cap D_i )
\overset{i_*}{\longrightarrow} H_1 ( \dot{X} )\oplus H_1( D_i )
\overset{j_*}{\longrightarrow} H_1(X_i)
\overset{\partial_*}{\longrightarrow}  \Tilde{H}_0( \dot{X} \cap D_i )
$
\end{enumerate}
Since $\dot{X} \cap D_i$ is connected, $\Tilde{H}_0( \dot{X} \cap D_i )=0$ and $j_*$ is a surjection.
Then, $H_1(X_i)
=\left( H_1 ( \dot{X} )\oplus H_1( D_i ) \right) / {\rm ker}j_*
=\left( H_1 ( \dot{X} ) \oplus 0 \right) / {\rm ker}j_*$.
Since the sequence is exact,
$H_1(X_i)=
\left( H_1 ( \dot{X} ) \oplus 0 \right) / {\rm im}i_*$.
We consider $\left<l_1\right>$, $\ldots$, $\left<l_n \right>$ as the part of a base of the first homology group $H_1( \dot{X} )$ of $\dot{X}$,
which is the free Abelian group generated by
$\left<l_1\right>$, $\ldots$, $\left<l_n \right>$,
$\left<p_1\right>$, $\ldots$, $\left<p_{r'} \right>$. ($r'=r(\dot{X})-n$.)
Since each sector is oriented, $[\partial D_i]= \sum_{k=1}^{n}d(l_k;e_i) \left< l_k \right>$ in $H_1( \dot{X})$ and ${\rm Im}i_*$ is generated by $\sum_{k=1}^{n}d(l_k;e_i) \left< l_k \right>$.
Then, 

\begin{enumerate}
\item [] $H_1(X_i)$
\item [] $=
\left(
\mathbb{Z} \left<l_1 \right>
\oplus \cdots
\oplus \mathbb{Z} \left< l_n \right>
\oplus \mathbb{Z} \left< p_1 \right>
\oplus \cdots
\oplus \mathbb{Z} \left< p_{r'} \right>
\right)
/ {\rm Span} \left\{ \sum_{k=1}^{n}d(l_k;e_i) \left< l_k \right> \right\}$
\item [] $=
\left(
\mathbb{Z} \left<l_1 \right>
\oplus \cdots
\oplus \mathbb{Z} \left< l_n \right>
\right)
/
{\rm Span} \left\{ \sum_{k=1}^{n}d(l_k;e_i) \left< l_k \right> \right\}
\oplus\\
\hspace{3mm}
\left(
\mathbb{Z} \left< p_1 \right>
\oplus \cdots
\oplus \mathbb{Z} \left< p_{r'} \right>
\right)$
\item [] $=
\left[
l_1, \ldots, l_n : \sum_{k=1}^{n}d(l_k;e_i)l_k
\right]
\oplus
\mathbb{Z}^{r'}$.
\end{enumerate}

Therefore,\\
$H_1(X)=
\left[
l_1, \ldots, l_n :
\sum_{k=1}^{n}d(l_k;e_1)l_k, \ldots, \sum_{k=1}^{n}d(l_k;e_m)l_k
\right]
\oplus \mathbb{Z}^{r'}$.
$\Box$
\end{pf}

\begin{figure}[htbp]
	\begin{center}
	\includegraphics[trim=0mm 0mm 0mm 0mm, width=.6\linewidth]{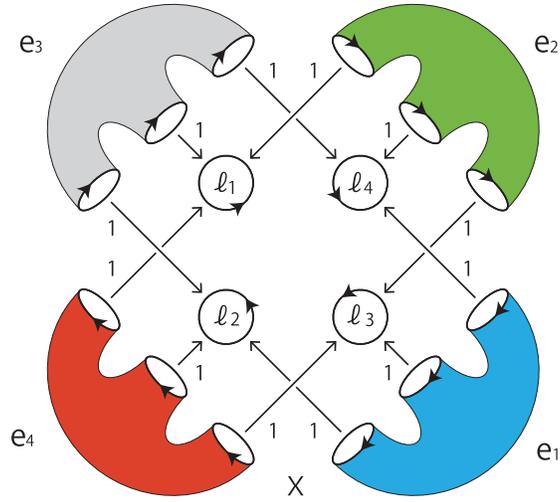}
	\end{center}
	\caption{A regular multibranched surface $X$ has four branches $l_1$, $l_2$, $l_3$ and $l_4$ and four sectors $e_1$, $e_2$, $e_3$ and $e_4$, which are $3$-punctured spheres. The oriented degree of every prebranch of $X$ is $1$.}
	\label{example1}
\end{figure}

\begin{Example}\label{mbs_example}{\rm 
For the multibranched surface $X$ in Figure \ref{example1},
$H_1(X)=\left( \mathbb{Z}/3\mathbb{Z} \right) \oplus \mathbb{Z}^4$
(The multibranched surface is defined in \cite{EMO}).
}
\end{Example}

\begin{pf}
The multibranched surface $\dot{X}$ is illustrated in Figure \ref{spine}.
By considering a spine of $\dot{X}$,
we have $H_1(\dot{X})=\mathbb{Z}^8$.
Since we can consider $\left< l_1 \right>$, $\left< l_2 \right>$, $\left< l_3 \right>$ and $\left< l_4 \right>$ as the part of generators of $\mathbb{Z}^8$,
we obtain the following equation  by Theorem \ref{thm_h1x}.
\begin{enumerate}
\item []
$H_1(X)$
\item []
$= \left[
l_1, l_2, l_3, l_4 :
     l_2+l_3+l_4,\
l_1     +l_3+l_4,\
l_1+l_2     +l_4,\
l_1+l_2+l_3
\right]
\oplus
\mathbb{Z}^4
$
\item []
$
= \left[
l_1, l_2, l_3, l_4 :
l_1,\
l_2,\
l_3,\
3l_4
\right]
\oplus
\mathbb{Z}^4
$
\item []
$
= \left[
l_4 :
3l_4
\right]
\oplus
\mathbb{Z}^4
=\mathbb{Z}/3\mathbb{Z} \oplus \mathbb{Z}^4
$.
\end{enumerate}

\begin{figure}[htbp]\begin{center}
\includegraphics[width=11cm]{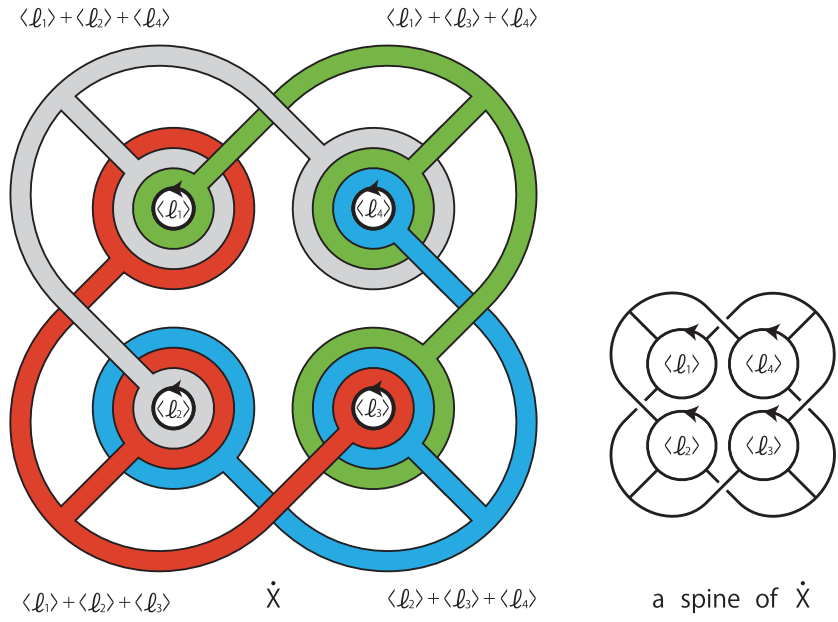}\end{center}
\caption{The multibranched surface $\dot{X}$ obtained from $X$ and a spine of $\dot{X}$.}
\label{spine}
\end{figure}

\end{pf}

\begin{Example}\label{1sector}{\rm
For the multibranched surface $X=L \cup_{\phi} E$ in Figure \ref{1sector},
$H_1(X)=\left( \mathbb{Z}/p \mathbb{Z} \right) \oplus \mathbb{Z}^{2g+n-1}$
, where
$L$ has $n$ components $l_1$, $\cdots$, $l_n$,
$E$ is a connected $n$-punctured closed surface with genus $g$, $\partial E = \{ c_1, \cdots, c_n \} $,
$p_i=d(c_i)$ and $p={\rm gcd}\{ p_1, \cdots, p_n\}$.
Moreover, if $p>1$, then $X$ cannot be embedded into $S^3$.
(cf. Example \ref{ex_gx2}.)
}\end{Example}

\begin{figure}[htbp]\begin{center}
\includegraphics[width=6cm]{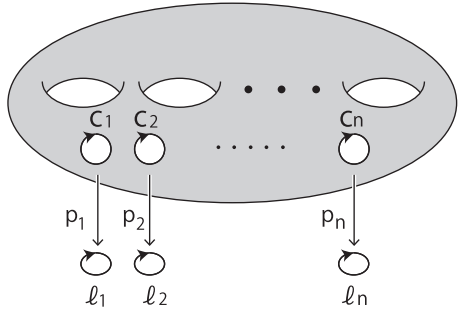}\end{center}
\caption{A multibranched surface obtained from one sector and $n$ branches.}
\label{1sector}
\end{figure}

\begin{pf}
In the same manner as Example \ref{mbs_example},
we obtain the following equation  by Theorem \ref{thm_h1x}.
\begin{enumerate}
\item []
$H_1(X)$
\item []
$= \left[
l_1, \cdots, l_n :
p_1  l_1 + \ldots + p_n  l_n
\right]
\oplus
\mathbb{Z}^{2g}
$
\item []
$
= \left[
l_1, \ldots, l_n :
pl_n
\right]
\oplus
\mathbb{Z}^{2g}
$
\item []
$
=
\left( \mathbb{Z}^{n-1} \oplus \left[ l_n : pl_n \right]  \right) \oplus \mathbb{Z}^{2g}
=\mathbb{Z}/p \mathbb{Z} \oplus \mathbb{Z}^{2g+n-1}
$ $\Box$
\end{enumerate}

\end{pf}

\noindent
{\bf Problem.}
When $p=1$, can a multibranched surface $X$ in Figure \ref{1sector} be embedded into $S^3$?

\begin{Remark}\label{embeddableintoS3}{\rm
When $p=1$ and $n \leq 2$,
a multibranched surface $X$ in Figure \ref{1sector}
can be embedded into $S^3$.
}\end{Remark}

\section{Minors of multibranched surfaces}\label{minors}
In Graph Theory, a finite graph $G'$ is called a {\it minor} of a finite graph $G$ if $G'$ is obtained from $G$ by removing and contracting some edges \cite{D}.
There are many researches with respect to minor of a graph, for example, Robertson and Seymour's Graph Minor Theorem has an important role \cite{L}.
We define minors for multibranched surfaces analogously and study properties.

\subsection{Intrinsically knotted (or linked) regular multibranched surfaces}\label{subsecIKIL}
A connected and closed $2$-dimensional manifold $F$ embedded into the $3$-dimensional sphere $S^3$ is {\it knotted} if $F$ is not a Heegaard surface of $S^3$. 
A disconnected and closed $2$-dimensional manifold $F$ embedded into $S^3$ is {\it linked} if there exists no essential sphere in $S^3 - F$.

\begin{Definition}\label{intrinsically knotted}{\rm
A regular multibranched surface $X$ embeddable in $S^3$ is {\it intrinsically knotted} (resp. {\it intrinsically linked}) if for the image of every embedding of $X$ into $S^3$,
it contains a knotted closed surface (resp. linked closed surfaces).
}\end{Definition}

\noindent A graph ($1$-dimensional CW complex) $G$ is {\it intrinsically knotted} (resp. {\it intrinsically linked}) if the image of every embedding of $G$ into $S^3$ contains non-trivial knot (resp.  non-trivial link).

\begin{figure}[h] 
\begin{center}
\includegraphics[trim=0mm 0mm 0mm 0mm, width=.9\linewidth]{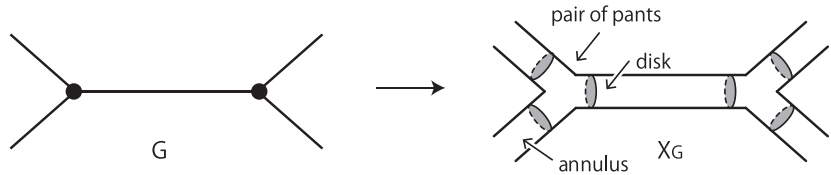}
\end{center}
\caption{A multibranched surface $X_G$ is obtained from $G$.}
\label{graph}
\end{figure}

\begin{Proposition}\label{propIK}
Let  $G$ be an intrinsically knotted (resp. intrinsically linked) graph.
Then the multibranched surface constructed as illustrated in Figure \ref{graph} is intrinsically knotted (resp. intrinsically linked).
\end{Proposition}

\begin{pf}
Let $E(G)$ (resp. $V(G)$) be a set of edges (resp. vertices) of $G$.
$X_G$ is obtained from $\#V$ spheres 
(corresponding to vertices of $G$) by attaching $\# E$ annuli (corresponding to edges of $G$).
Suppose that $X_G$ is embedded into $S^3$.
Since there exists a handlebody $W$ such that
$\partial W \subset X_G$ and
a spine of $W$ is homeomorphic to $G \subset S^3$,
a torus in $X_G$ corresponding to $k$ is knotted in $S^3$, where $k$ is a non-trivial knot contained in $G$.
Therefore $X_G$ is intrinsically knotted (linked). 
$\Box$.
\end{pf}

\begin{figure}[htbp]\begin{center}
\includegraphics[trim=0mm 0mm 0mm 0mm, width=.6\linewidth]{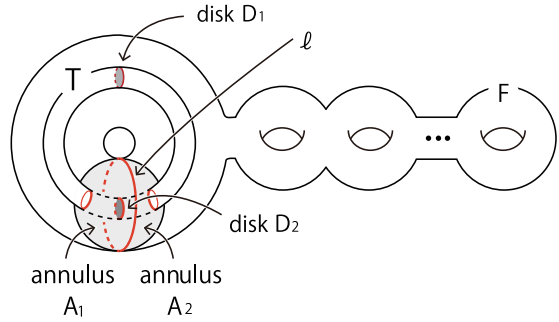}\end{center}
\caption{Intrinsically linked surface, in which every branch has index $1$.}
\label{linked}
\end{figure}

\begin{Example}\label{ex_IL}{\rm
The following multibranched surface $X$ in Figure \ref{linked} is intrinsically linked.
}\end{Example}

\begin{pf}
For every embedding of $X$,
there exists a solid torus $V$ such that $\partial V \subset T$, $D_1 \subset V$ and $D_2 \subset V$.
Since $l \subset F$ is contained in $S^3 - V$ and $l$ is parallel to a meridian of $V$, there is no sphere which separates $T$ from $l$.
Therefore the closed surfaces $T \cup F$ is linked.
$\Box$\end{pf}

\subsection{Topological minors of multibranched surfaces}

Let $X=L \cup_\phi E$ be a multibranched surface and let $e$ be a sector of $X$.
We define a continuous map $\phi' : \partial (E-e) \to L$ as follows:
for every connected boundary $c$ of $E-e$,
$\phi'|_c=\phi|_c$.
The multibranched surface $X'=L \cup_{\phi'} (E-e)$ obtained from the triple $(L, E-e, \phi')$ is said to be obtained by {\it removing} a sector $e$ of $X$. 

Let $A = S^1 \times [-1, 1]$ be an annulus of a multibranched surface $X=L \cup_\phi E$ with $A \cap L = \partial A$.
(Putting $l^- = S^1 \times \{-1 \}$ and $l^+ = S^1 \times \{ 1 \}$.)
By identifying an interval $\{ s \}  \times [-1, 1] \subset A$ with a point $\{ s \} \times \{ 0 \} \subset A$,
we obtain a multibranched surface,
which is said to be obtained by {\it contracting} an annulus $A$ of $X$.

Let $X$ and $Y$ be multibranched surface.
If there exist multibranched surfaces $X'$ homeomorphic to $X$ and
$Y'$ homeomorphic to $Y$ such that
$X'$ is obtained by removing a sector (resp. by contracting an annulus) of $Y'$,
we denote $X\overset{\mathrm{r}}{<} Y$ (resp. $X\overset{\mathrm{c}}{<} Y$).
Also, we denote $X < Y$ if
$X\overset{\mathrm{r}}{<} Y$ or $X\overset{\mathrm{c}}{<} Y$.

The set consisting of all multibranched surfaces is denoted by $\mathcal{M}$.
We define an equivalence relation $\sim$ on $\mathcal{M}$ as follows:
if $X < Y$ and $Y < X$, then $X \sim Y$.
An element of the quotient set $\mathcal{M}/\sim$ is called a multibranched surface {\it class} (or a multibranched surface for simplicity).
We define a partial order $\prec$ over $\mathcal{M}/\sim$ as follows.

\begin{Definition}\label{def_top_minor}{\rm
Let $X$ and $Y$  be multibranched surfaces.
We denote $[X] \prec [Y]$ if
there exists a finite sequence $\{ X_0, X_1, \ldots, X_{n-1}, X_n \}$ of multibranched surfaces such that
$X_0 \sim X$, $X_n \sim Y$ and $X_0 < X_1 < \cdots < X_{n-1} < X_n$.
}\end{Definition}

A multibranched surface (class) $[X]$ is said to be a {\it minor} of a multibranched surface (class) $[Y]$
if $[X] \prec [Y]$.
Also, $[X]$ is said to be a {\it proper minor} of $[Y]$
if $[X] \prec [Y]$ and $[Y] \not= [X]$.
A subset $\mathcal{P}$ of $\mathcal{M}/\sim$ is {\it minor closed}
if for every multibranched surface $X \in \mathcal{M}$, every minor of $[X]$ belongs to $\mathcal{M}/\sim$.
For minor closed set $\mathcal{P}$, we define the {\it obstruction set} $\Omega(\mathcal{P})$ as follows.
\[
\Omega(\mathcal{P})=
\left\{
[X] \in \mathcal{M}/\sim | [X] \not \in \mathcal{P}, \mbox{ Every proper minor of } [X]  \mbox{ belongs to } \mathcal{P}. 
\right\}
\]

\begin{Proposition}\label{minorclosed}
The following sets are minor closed.
\begin{enumerate}
\item The set consisting of multibranched surfaces embeddable into $S^3$,
denoted by $\mathcal{P} {S^3}$
\item The set consisting of multibranched surfaces embeddable into $S^3$ which is not intrinsically knotted (resp. linked),
denoted by $\mathcal{P}_{knot}$ (resp. $\mathcal{P}_{link}$).
\end{enumerate}
\end{Proposition}

\begin{pf}
1.
Let $X$ be a multibranched surface embeddable into $S^3$ and
let $Y$ be a multibranched surface with $Y < X$.\\
(Case 1: $Y \overset{\mathrm{r}}{<} X$.)
Since $X$ is embeddable into $S^3$, the subspace $Y$ of $X$ is embeddable into $S^3$.\\
(Case 2: $Y \overset{\mathrm{c}}{<} X$.)
Suppose that $X$ is embedded in $S^3$.
We can contract $A$ to the core of $A$ in $S^3$.
Therefore $Y$ is embeddable into $S^3$.

2.
Let $X$ be a non intrinsically knotted (resp. linked) multibranched surface and
let $Y$ be a multibranched surface with $Y < X$.
Suppose that $X$ is embedded into $S^3$ and $X$ has no knotted (resp. linked) surface.\\
(Case 1: $Y \overset{\mathrm{r}}{<} X$.)
If $Y$ has knotted (resp. linked) surface, $X$ has knotted (resp. linked) surface. This contradicts our assumption.\\
(Case 2: $Y \overset{\mathrm{c}}{<} X$.)
By contracting $A$ to the core of $A$ in $S^3$,
we obtain $Y \subset S^3$ which has no knotted (resp. linked) surfaces.
$\Box$.
\end{pf}

In \cite{EMO}, we give some examples of multibranched surfaces which belong to $\Omega(\mathcal{P} {S^3})$.

\begin{Example}\label{ex_critical}{\rm
Let $X$ be a regular multibranched surface in Figure \ref{1sector}.
If $p={\rm gcd} \{p_1, p_2, \ldots, p_n \}$ is not $1$, then
$[X] \in \Omega(\mathcal{P} S^3)$.
}\end{Example}

\begin{pf}
By \cite{EMO}, if $p$ is not $1$, $X$ is not embeddable into $S^3$.
Since $\dot{X}$ is embeddable into $S^3$,
$[X] \in \Omega(\mathcal{P} S^3)$
$\Box$
\end{pf}

\begin{figure}[htbp]\begin{center}\includegraphics[trim=0mm 0mm 0mm 0mm, width=.5\linewidth]
{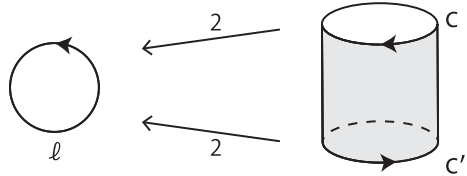}\end{center}
\caption{$g(X)=2$. $X$ has one branch and one sector and two prebranches $c$ and $c'$ with $od(c)=od(c')=2$.}
\label{obstruction}
\end{figure}

\begin{Example}\label{ex_gx2}{\rm
For the following multibranched surface $X$ in Figure \ref{obstruction}, $[X] \in \Omega(\mathcal{P}{S^3})$ and $g(X)=2$.
}\end{Example}

\begin{pf}
Suppose that $X$ is embeddable into a lens space $L(p,q)$.
For every neighborhood $N \in \mathcal{N}(X)$ of $X$,
$N$ is a Seifert manifold with a M\"{o}bius band as its base and a single singular fiber.
Let $L(p,q)=N \cup V$ ($N \cap V=\partial N = \partial V$).
We note that the neighborhood $\partial N$ is compressible in $L(p,q)$.
Since a lens space has no essential and separating sphere,
$V$ is a solid torus.
Then, $L(p, q)$ is a Seifert manifold with a projective plane $\mathbb{R}P^2$ as its base and $1$ or $2$ singular fibers.
If there are $2$ singular fibers,
there exists an incompressible torus.
This leads to a contradiction.
Therefore there is a $1$ singular fiber.
By [Proposition 10.11, \cite{FM}], $L(p,q)$ is homeomorphic either to a Seifert manifold with a sphere as its base and $3$ singular fibers or 
homeomorphic to $\mathbb{R}P^3 \# \mathbb{R}P^3$.
This leads to a contradiction.
Therefore, $X$ is not embeddable into any lens space.
One can show that $X$ is not embeddable into $S^2 \times S^1$.
$\Box$
\end{pf}

\begin{Example}\label{ex_link}{\rm
For the multibranched surface $X$ in Figure \ref{linked},
let $X'$ be a multibranched surface obtained by removing open disks from $A_1$ and $A_2$. Then, $[X'] \in \Omega(\mathcal{P}_{link})$.
}\end{Example}

\begin{pf}
Since for every embedding of $X'$ into $S^3$,
there is a loop $k$ contained in $T$ such that the linking number of the link $k \cup l$ is not $0$, $F \cup T$ is linked surface.
Since every multibranched surface obtained by removing a sector
(or contracting an annulus) contains no linked surface or it is equivalent to $X'$,
$[X'] \in \Omega(\mathcal{P}_{link})$.
\end{pf}

\subsection{Neighborhood minors of regular multibranched surfaces}

Let $X=L \cup_\phi E$ be a multibranched surface and let $l$ be a branch of $X$.
We define a continuous map $\phi' : \partial E \to L$ as follows:
For every connected boundary $c$ of $E$ attached to $l$,
(i) $\phi'|_c : c \to l$ is a covering map with covering degree $1$ and (ii) $\phi'|_c$ and $\phi|_c$ have same orientation.
For every connected boundary $c$ not attached to $l$, $\phi'|_c = \phi|_c$.
The multibranched surface $X'=L \cup_{\phi'} E$ obtained from the triple $(L, E, \phi')$ is said to be obtained by {\it reducing a degree} of a branch $l$ of $X$. 
We define a binary relation $\precN $ on $\mathcal{M}$ as follows.

\begin{Definition}\label{minor2}{\rm
For multibranched surfaces $X$ and $Y$, we denote $X \precN Y$ if
for every neighborhood $N \in \mathcal{N}(Y)$,  $X$ is embeddable into $N$.
This relation is a preorder on $\mathcal{M}$. If $X \precN Y$, then $X$ is said to be a {\it neighborhood minor} of $Y$.
}\end{Definition}

\begin{Remark}\label{rem_nbd_minor}{\rm
For regular multibranched surfaces $X$ and $Y$,
if $X \precN Y$, then $g(X) \leq g(Y)$.
}\end{Remark}

\begin{Proposition}\label{minorclosed}
If a regular multibranched surfaces $X$ and $Y$
are satisfying one of the following conditions, then $X \precN Y$.
\begin{enumerate}
\item $X$ is obtained by contracting an annulus of $Y$.
\item $X$ is obtained by reducing a degree of a branch of $Y$.
\item $X$ is obtained by a connected sum of a torus to a sector of $Y$.
\end{enumerate}
\end{Proposition}

\begin{pf}
Let $N$ be a neighborhood of $Y$.

1.
Suppose that $A$ be an annulus of $Y$ and $X$ is obtained by contracting $A$ of $Y$. 
Let $l^*$ be a branch of $X$, which occurs by contracting $A$.
$X$ is embeddable into $N$ such that $l^*$ is contained in a solid torus $A \times [-1, 1] \subset N$, which is corresponding to $A$ of $Y$.
See Figure \ref{contraction}.

\begin{figure}[htbp]
	\begin{center}
	\includegraphics[width=8cm]{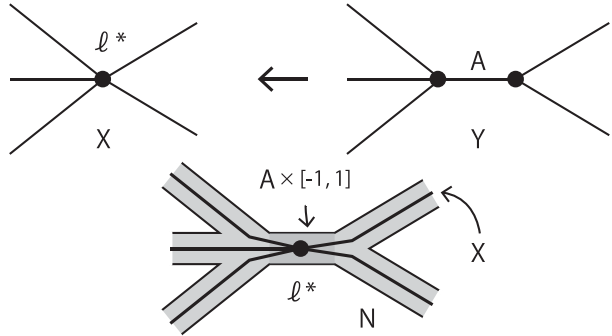}
	\end{center}
	\caption{$X$ is embeddable into a neighborhood $N$ of $Y$.}
	\label{contraction}
\end{figure}

2.
Since $Y$ is embeddable into $N$ (such that $Y$ is a spine of $N$),
let $Y$ be a multibranched surface embedded into $N$.
Let $l$ be a branch of $Y$.
We reduce a degree of $l$ of $Y$ into $1$ as follows.
Put $T=\partial N(l)$.
By cutting $T$ along the sector of $Y$,
we obtain $n$ annuli $A_1, \ldots, A_n$, where $n=i(l)$.
Let $C_1, \ldots, C_n$ be the collar of $E$ which is incident to $l$.
We modify $C_j$ to $A_1 \cup \cdots \cup A_{j-1}$ ($j=2, \ldots, n$)
so that they are mutually interior-disjoint and whose boundary is $l^*=C_1 \cap T$.
See Figure \ref{reducing}.
Then, in the modified multibranched surface $X$,
we have $d(l^*)=1$.

\begin{figure}[htbp]
	\begin{center}
	\includegraphics[width=8cm]{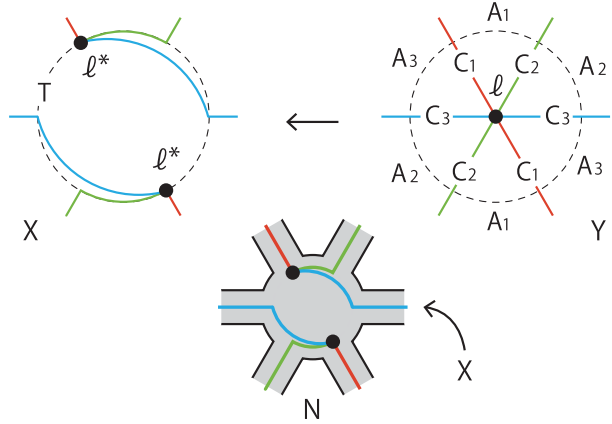}
	\end{center}
	\caption{In the case where $d(l)=2$ and $i(l)=3$. }
	\label{reducing}
\end{figure}

3.
Let $Y$ be a multibranched surface embedded into $N$.
By a connected sum of a ``small'' torus to a sector of $Y$,
we obtain an embedding of $X$ into $N$.
$\Box$\end{pf}


\end{document}